# Mathematical Foundations of Consciousness

Willard L. Miranker[1], Gregg J. Zuckerman[2]
Departments of [1]Computer Science, [2]Mathematics
Yale University

5/15/08

This work is dedicated to the memory of our mentors.

*"Mathematics as an expression of the human mind reflects the active will, the contemplative reason, and the desire for aesthetic perfection. Its basic elements are logic and intuition, analysis and construction, generality and individuality. Though different traditions may emphasize different aspects, it is only the interplay of these antithetic forces and the struggle for their synthesis that constitute the life, usefulness, and supreme value of mathematical science."* [1]*Richard Courant* (1941).

*"One expects that logic, as a branch of applied mathematics, will not only use existing tools from mathematics, but also that it will lead to the creation of new mathematical tools, tools that arise out of the need to model some real world phenomena not adequately modeled by previously known mathematical structures."* [2]*Jon Barwise* (1988).

**Abstract:** We employ the Zermelo-Fränkel Axioms that characterize sets as mathematical primitives. The Anti-foundation Axiom plays a significant role in our development, since among other of its features, its replacement for the Axiom of Foundation in the Zermelo-Fränkel Axioms motivates Platonic interpretations. These interpretations also depend on such allied notions for sets as pictures, graphs, decorations, labelings and various mappings that we use. A syntax and semantics of operators acting on sets is developed. Such features enable construction of a theory of non-well-founded sets that we use to frame mathematical foundations of consciousness. To do this we introduce a supplementary axiomatic system that characterizes experience and consciousness as primitives. The new axioms proceed through characterization of so-called consciousness operators. The Russell operator plays a central role and is shown to be one example of a consciousness operator. Neural networks supply striking examples of non-well-founded graphs the decorations of which generate associated sets, each with a Platonic aspect. Employing our foundations, we show how the supervening of consciousness on its neural correlates in the brain enables the framing of a theory of consciousness by applying appropriate consciousness operators to the generated sets in question.

## 1. Introduction

Analytic writing on mind and consciousness dates to Aristotle's De Anima (Ross, ed. 1961). Yet to this day the phenomena of consciousness continue to elude illuminating scientific characterization. We should not be surprised at this since,

> *"A physical scientist does not introduce awareness (sensation or perception) into his theories, and having thus removed the mind from nature, he cannot expect to find it there."* Schrödinger 1958

The self-referential qualities of consciousness place it outside conventional logic(s) upon which scientific models and frameworks have heretofore been constructed. However more contemporary mathematical development has begun to deal with features of self-reference. We shall address Schrödinger's critique by assembling and extending such development thereby putting self-reference as a form of awareness into theory. In this way we shall frame mathematical foundations for a theory of consciousness. Then as an application to a neural network model of brain circuitry, we shall exhibit a theory of consciousness using these foundations.

### 1.1 Mathematical thought and its limits
Platonism, that is, the interplay of ideal[1] and physical worlds, characterizes a central feature of mathematical thought. The briefest summary of the evolution of this Platonic dualism in mathematical thought and modeling might be made by citing the contributions of Euclid (the axiomatic method), Aristotle (the law of the excluded middle), Cantor 1895 (set theory), Russell ~1899 (his well-known paradox in set theory), Zermelo 1908 and Fränkel 1922 (the axioms of set theory that serve to accommodate Russell's paradox) and Gödel (incompleteness, a self-referential development). We shall extend this line and employ the axiomatic theory of sets to further characterize self-referential features.

The work of Zermelo-Fränkel and others transformed sets from so-called naïve objects into mathematical primitives (i.e., ideal Platonic objects). The Russell paradox and its accommodation demonstrate limitations on mathematical thought (about sets and related constructs). Today we are not surprised by such a limitation, since we have the well-known example of Heisenberg. The Platonic character of the latter is characterized by the Heisenberg inequality, its ideal form (Dym, McKean 1972), and its real world character by the limitation on the accuracies with which certain concurrent measurements can be made. The quality of self-reference (set self-membership) underlying the Russell paradox informs development of the ideal Platonic structures (i.e., of placing awareness into theory) required for constructing the mathematical foundations we seek.

### 1.2 Consciousness and its limits
As with the self-referential potentiality of naïve set theory, the self-referential character of consciousness appears paradoxical. It seems to be an illusion. The incompleteness of mathematical thought demonstrated by Gödel, suggests that all thought, and so consciousness in particular, is not explainable via a conventional approach such as by a Turing machine computation (Penrose 1989). *Incompleteness, while precluding establishment of certain knowledge within a system, allows for its establishment by looking onto the system from the outside. This knowledge from the outside (a kind of observing) is reminiscent of consciousness that provides as it does a viewing or experiencing of what's going on in thought processing.* Note the correspondence of these observations to Freud's meta-psychology where he recognizes a disconnect between mental and physical states,

> "...*mental and physical states represent two different aspects of reality,
> each irreducible to the other* (Solms 1994)."

---

[1] An ideal object or concept in the Platonic sense will sometimes be referred to as a Platonic object or concept or for emphasis as an ideal Platonic object.

However we may say that Freud's psychoanalytic method is a tool devised for penetrating the mental from the outside via the physical. Compare Freud's dual aspects of reality with the Platonic pairing of Descartes 1637, namely the *res cogitans* (ideal) and the *res extensa* (physical).

To frame a set theoretic correspondent to these features note that in axiomatic theory, a set has an inside (its elements) and an outside (the latter is not a set, as we shall see), and this allows a set to be studied from the outside. We liken this to interplay between the ideal (Platonic) and physical (computable) worlds, the latter characterizing a model for study from the outside of the former. So we expect consciousness to be accessible to study through extensions of the self-reference quality characterized by axiomatic set theory, in particular, by a special capacity to study a set from the outside. We do not claim that this gives a complete characterization of consciousness, although it might very well do so in the end. Rather this approach is an effectual way to introduce awareness into a theory (accommodating thereby Schrödinger's critique) and so to penetrate this elusive phenomenon.

**1.3 Summary**

Sect. 2 begins with a description of the crises in mathematical thought precipitated by Cantor's set theory and characterized by the Russell paradox. We describe how Gödel's discoveries inform the crises and furnish motivation for our development. We introduce a mathematical framework that includes sets, graphs, decorations, and the notion of non-well-founded sets and which enables annunciation of the anti-foundation axiom of set theory. This axiom allows replacement of the Russell paradox by a logically coherent dichotomy and is key to framing our approach characterized by observation of sets from the outside.

In Sect. 3 we introduce the Russell operator $\mathcal{R}$, a distinguisher between so-called normal and abnormal sets. A number of properties of $\mathcal{R}$ is collected, these to play a central role in the foundations to be developed. Then we introduce a number of other operator mappings along with interrelations, these to supplement $\mathcal{R}$ in the analysis of sets to follow. This operator syntactic framework is followed by a semantic development in which experience and consciousness are introduced as primitives. A Semantic Thesis for consciousness is then proposed, and a list of axioms for associated operators along with a descriptive semantics for each axiom is given (compare Aleksander, Dunmall 2003). The axioms along with their semantics are used for characterizing both the primitives and consciousness. $\mathcal{R}$ is shown to satisfy the axioms, giving it thereby the role of a so-called consciousness operator. This existence of a consciousness operator establishes consistency of the new axioms. Examples both of sets and operators illustrating the syntax and semantics are given.

   In Sect. 4 we give a description of tools for building a theory of consciousness upon the foundations developed. This begins with a formal process for labeling and then decorating a graph. The process establishes a way to induce the existence of a virtual set associated intrinsically with a graph (a two-level or self-referential feature). We then introduce a mapping construct called a histogram, a tool for applying this set with graph association process to a special class of graphs arising in brain circuitry. The M-Z equation is then developed, this equation characterizing a method for specifying the intrinsic set in question, including those that arise in brain circuitry. Finally the theory of consciousness is formulated as an application in which we employ neural network theory (Hebb's rule for synaptic weight change and the McCulloch-Pitts equation for neuronal input-output dynamics, (see Haykin 1999)) to specify the special class of labeled graphs in question. This two-level procedure is interpreted as a Platonic process (that is, the association of a virtual set with a graph) by means of what we call a Neural Network Semantic Thesis. An example of a neural state that instantiates the concept of a particular natural number is given. To complete the description of information processing from sensory perception through to consciousness, a third, purely physical, so-called Neuro-physiological Thesis is introduced. Sect. 4 concludes with a critical description both of these three theses and the analytic formalisms developed earlier. This critique serves to illuminate the mathematical foundations of consciousness developed.

   In Sect. 5 we ascribe syntactic and semantic nomenclature to a collection of basic operators, also offering interpretations of the role each plays in our theory. The flow of information from sensory input to conscious experience is described. Speculation is offered on the role of the sets we have constructed in this information flow. Finally a class of operators that characterize qualia is described.

   In Sect. 6 directions of future work are laid out. These include (i) examples and applications of the M-Z equation developed in Sect. 4, including the development of associated dynamics induced by the consciousness operators, (ii) the study of the diagonalization of $\mathcal{K}_A$, a special consciousness operator that informs the study of qualia and their neural correlates, (iii) connection of our mathematical foundations with processes of evolution, (iv) study of bi-simulation of graphs that characterizes the case that two memes share a thema, (v) model theoretic foundations of Aczel theory dealing with the consistency of the Z-F axioms with anti-foundation replacing foundation, (vi) study of the algebra of the fundamental operators appearing in Table 3.1, and (vii) classification of the consciousness operators and the connection of doing this to Gödelian incompleteness.
   The axioms of set theory that we employ explicitly are given in an appendix. This is followed by a glossary.

## 2. Preliminaries

In this section we describe the crises in mathematical thought engendered by the notions developed by Cantor, Russell and others. Then we describe the evolution of the crises according to the development of Zermelo-Fränkel, Gödel and others. We continue with the introduction of terminology and properties that provide the setting for our work.

### 2.1 Crises in mathematical thought
We begin with Cantor's definition of a set, often regarded as the *naïve notion of set*.

> "*A set is a collection into a whole of definite, distinct objects of our intuition or <u>thought</u>.*"

When specificity is required, we shall hereafter use the term *collection* for a set in the sense of Cantor's definition.

Cantor's use of the word "thought" shows that set theory is entwined with consciousness from the start. In fact, Cantor's definition is circular, replacing one mystery by another. It replaces the unanswered questions: what is a definite object? what is thought? by others, namely: who does the collecting? the thinking? The latter have a correspondence to the questions often raised in consciousness studies, "Who is doing the looking? the experiencing?" Suppose the words "intuition or thought" in Cantor's definition are replaced by the word "consciousness". This would make it an exception to Schrödinger's critique, relating it to what is perhaps the only other known exception, namely to Von Neumann's (mysterious) appeal to the observer's consciousness of the outcome of a measurement to specify the moment of collapse of the wave function during a quantum mechanical measuring process.

Cantor's definition of a set supports a logical inconsistency, resulting in several paradoxes. The most accessible of these is the Russell paradox that goes to the essence of that inconsistency. This paradox is expressed in terms of the *Russell naïve set N*, which is the collection of all sets $x$ such that $x$ is not a member of $x$. The logical inconsistency of $N$ is revealed by the following observations:

1. Since $N$ is a set, either $N \in N$ or $N \notin N$.
2. If $N \in N$, then $N \notin N$. If $N \notin N$, then $N \in N$. (2.1)

The annunciation of this paradox by Russell (Zermelo is thought to have known earlier of the paradox) precipitated a major crisis in mathematical and philosophical thought. In 1893, Frege had just completed development of an axiomatic treatment of sets when a letter to him from Russell informing him of the paradox overturned his central thesis. Various mathematicians (Bernays, Gödel, Hilbert, Russell, Von Neumann, Whitehead…) attempted to rework the foundations of mathematics so as to resolve the(se) paradox(es). It is the axiomatic approach to set theory that provides for us the most fruitful resolution, motivating our own development. (See the appendix for these axioms.) The key feature of the axiomatic approach is to regard the concept "set" as a primitive (an undefined notion), and the concept "is an element of" as a primitive relation. The axioms are chosen to ensure that there does not exist a set $y$ such that $x \in y$ if and only if $x \notin x$; in other words, *within axiomatic set theory, there is no Russell set.* Even so, this axiomatic approach allows for a coherent elaboration of the quality of self-reference in set theory, and so, it supports the connection of the study of sets to the development of the mathematical foundations we are after.

We use Z-F, the Zermelo-Fränkel axioms of set theory, however replacing FA, the Foundation Axiom (a latter day addition by Von Neumann to the original Z-F list) by AFA, the Anti-foundation Axiom. To distinguish a set in the sense of these axioms from a collection of Cantor, we shall use the terminology, *bona fide set* for the former.

Although successfully accommodating the paradox, the axiomatic development of set theory brought with it a deeper problem: is the axiomatic system itself consistent? That is, can we derive a logical inconsistency from the axioms? Gödel produced a two level approach to this issue. At a mathematical level is a set theoretic formula, and at a meta-mathematical level is the proposition asserting the consistency of set theory. We interpret this as an instance of self-reference, a viewing of a mathematical object meta-mathematically, that is from the outside. Gödel showed that if axiomatic set theory is consistent then it is incomplete. This incompleteness is widely celebrated (see Gödel-Escher-Bach of Hofstadter 1979, Emperor's New Mind of Penrose 1989, Scientific American 1968…).

One might say that Gödel replaced one crisis in mathematical thought by another. Subsequently, mathematicians (Aczel, 1988…) did show that if Z-F with FA deleted is consistent, then Z-F with AFA replacing FA is also consistent. These results of Gödel and his successors provide for us the framework to develop our self-referential two level approach that consists, in particular, of *a syntactic level and a semantic level.*

**2.2 Sets, graphs, decorations, the axiom of anti-foundation**
The special nature of set theory can be traced in part to the use of two different notions of belonging associated with sets. One is denoted by $\in$ (the primitive concept '*is an element of*') and the other by $\subset$ (for the concept '*is a subset of*').

For clarity we adopt the following notational conventions.

a) Sets will be denoted by Latin characters, *a, A, b*…
Braces will also denote a set, the contents of which and/or conditions specifying the set placed within the braces: $\{\textit{list of set elements and/or conditions for being a set element}\}$.
b) Mappings between sets will be denoted by lower case Greek characters, $\alpha, \beta$...
c) Relations and operators as well as certain special objects to be introduced called classes will be denoted with upper case script Latin letters, $\mathcal{A}, \mathcal{B},…\mathcal{R}…$ A generic operator will be denoted by an upper case script $O$.
d) The empty set $\{x \mid x \neq x\}$ will, as usual, be denoted by $\varnothing$. The existence of $\varnothing$ follows from the Z-F Axioms of Existence and Comprehension (see the appendix).

We shall restrict our attention to pure sets specified as follows.

**Definition 2.1:** A set is a *pure set* if its elements are sets, the elements of its elements are sets…

Note that any finite collection (naïve set) of objects that are not themselves bona fide sets furnishes an example of a not pure set. Our presentation involves both normal and abnormal sets, these set types specified as follows.

**Definition 2.2:** A set *x* is *normal* if $x \notin x$. It is *abnormal* if $x \in x$.

The Quine atom specified in the following definition supplies and example of an abnormal set.

**Definition 2.3:** The *Quine atom* $\Omega$ is the set defined by the condition $\Omega = \{\Omega\}$.

We shall make use of a collection of notions specified in the following paragraph. (See Aczel 1988, Chap. 1.)

A *graph*[2] will consist of a collection *N* of *nodes* and a collection *E* of *edges,* each edge being an ordered pair $(n, n')$ of nodes. No knowledge of the nature of the elements of *N* is required. If $(n, n')$ is an edge, we shall write $n \to n'$ and say that $n'$ is a *child* of the node *n,* the latter called a parent of the node $n'$. A *path* is a sequence (finite or infinite)

$$n_0 \to n_1 \to n_2 \to \cdots$$

of nodes $n_0, n_1, n_2 \ldots$ linked by edges $(n_0, n_1), (n_1, n_2) \ldots$ A *pointed graph* is a graph together with a distinguished node called its *point*. A pointed graph is accessible, i.e., is an *accessible pointed graph* (apg), if for every node *n* there is a path $n_0 \to n_1 \to \cdots \to n$ from the point $n_0$ to the node *n*. If this path is always unique then the pointed graph is a *tree,* and the point is the *root* of the tree. A *decoration* of a graph is an assignment of a set to each node of the graph so that the elements of the set assigned to a node are the sets assigned to the children of that node. Alternatively a decoration is a set valued function *d* on *N* such that

$$\forall a \in N, \quad da = \{db \mid a \to b\}. \tag{2.2}$$

A *picture* of a set is an accessible pointed graph that has a decoration in which the set is assigned to the point. A given set may have many pictures. Being well-founded, a key property of graphs and sets is specified in the following definition.

**Definition 2.4:** A graph is *well-founded* if it has no infinite path. It is *non-well-founded* otherwise..

---

[2] What we call a graph is in fact a directed graph. For convenience we drop the descriptor directed throughout.

With this terminology, we collect the known results stated in the following proposition.

**Proposition 2.5:**  *i*) Every well-founded graph has a unique decoration.
*ii*) Every well-founded apg is a picture of a unique set.
*iii*) Every set has a picture.

Continuing, we define well-foundedness for sets.

**Definition 2.6:**  A set is *well-founded* if its picture is well-founded. It is *non-well-founded* otherwise.

An alternate name for a non-well-founded set is a hyper-set, but we prefer never to use the latter term.

We now state the *anti-foundation axiom* that is central to the development. Note it is stated for general graphs that are neither necessarily pointed nor necessarily accessible.

**AFA** (Aczel)**:** *Every graph has a unique decoration.*

Some consequences of this axiom are given in the following proposition.

**Proposition 2.7:** 1. Every pointed graph is the picture of a unique set.
2. Non-well-founded sets exist.
3. A non-well-founded graph will picture a non-well-founded set.
4. Every set is the decoration of at least one apg.

**Proof:**  See Aczel 1988.

The relationship between these concepts is summarized in terms of two mappings, the tree mapping $\tau$ and the decoration of the point $P$ mapping $\delta$ is shown in Fig. 2.1.

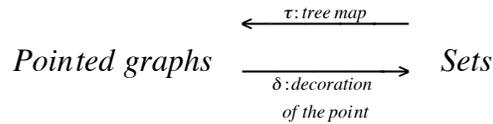

**Figure 2.1:** Schematic of the mappings $\tau$ and $\delta$

There are many graphs $\Gamma_i$, the decoration of whose point is a given set $A$. That is, for the map $\delta$, we have

$$\delta\Gamma_1 = \delta\Gamma_2 = \cdots = A. \tag{2.3}$$

However there is a unique pointed graph, $\Gamma_* = \Gamma_*(A)$ called the canonical tree of $A$, such that $\delta\Gamma_* = A$ and

$$\tau A = \Gamma_*(A). \tag{2.4}$$

The canonical tree of a set is specified in the following definition.

**Definition 2.8:** A finite collection $x, x_1, \ldots, x_n$ of sets forms a chain beginning at $x$ if $x_n \in x_{n-1} \in \cdots \in x_1 \in x$. The tree $\tau x$ of $x$ is the graph whose nodes are chains beginning at $x$ and whose edges are given by $(x_n \in x_{n-1} \in \cdots \in x_1 \in x, \; x_{n+1} \in x_n \in \cdots \in x_1 \in x)$.

$\delta(\Gamma, p)$ will denote the set associated with the node $p$ of the graph $\Gamma$ in the decoration of the latter. So $\delta\Gamma = \delta(\Gamma, P)$ is the set in the decoration of the pointed graph $\Gamma$ that corresponds to the point $P$ of $\Gamma$. Then a sufficient condition for normality of a set is given in the following proposition.

**Proposition[3] 2.9:** If for every child $c$ of $P$, $\delta(\Gamma, c) \neq \delta(\Gamma, P)$, then $\delta(\Gamma, P)$ is normal.

### 2.3 Classes and mappings

Classes are primitives introduced by Gödel. A collection of sets with a common property is called a *class*. A set is also a class; a class that is not a set is called a *proper class*. The elements of a class are sets, the sets being the primitives defined by the Z-F axioms with the AFA replacing the FA. Conversely, any set is a member of a set.

We now formalize the notions of several types of classes to be used. These are: relations, functions, and operators. They are illustrated by the nest of concepts shown in Fig. 2.1, the outermost member of which is comprised of the classes.

Inside of classes is the collection of relations. A *relation* is a class consisting of ordered pairs of sets.

Inside of relations is the collection of functions. A *function* is a relation with the graph property: namely, if $(x, y)$ and $(x, z)$, both being ordered pairs of sets in a relation $\mathcal{F}$, implies that $y = z$, then $\mathcal{F}$ is said to have the *graph property*.

Inside of functions is the collection of operators. An *operator* $O$ is a function whose domain is the class of all sets. An operator is a relation, since $O x$ is the unique $y$ so that $(x, y) \in O$.

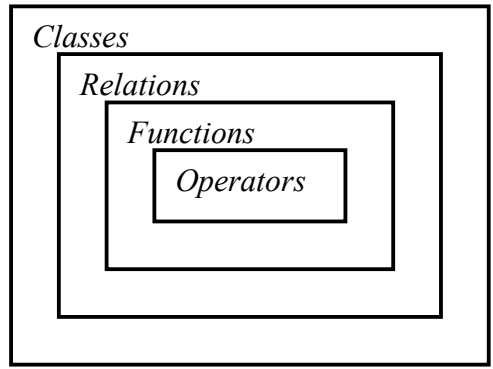

**Figure 2.1:** Nesting within the theory of classes

---

[3] Proofs that follow directly from definitions are omitted throughout.

## 3. The Russell Operator, Operator Syntax and Semantics, Semantic Thesis, Axioms

In this section we supply syntax and semantics for some operators of relevance for our axiomatic treatment of consciousness. We start in Sect. 3.1 with the characterization of the Russell operator, since it plays a central role. Then in Sect. 3.2, we introduce a relevant collection of operators and develop mathematical properties (syntax) for them. In Sect. 3.3, we state the Semantic Thesis that characterizes consciousness as the action of operators on experience. *Consciousness and experience are introduced as primitives*, and an open axiom system for them is elaborated. The axioms are accompanied by semantic characterizations of the associated operators.

### 3.1 The Russell operator

The Russell operator $\mathcal{R}$ plays a special role in the syntax and semantics of the development of the Semantic Thesis. $\mathcal{R}$ is defined by its action on a set $A$ as follows.

**Definition 3.1:** $\mathcal{R} A = \{x \in A \mid x \notin x\}$.

We see that $\mathcal{R}$ may be viewed as a selector of the normal elements of $A$ and a rejecter of the abnormal. A formal definition of a selector (operator) will be given in Def. 3.14. $\mathcal{R}$ is a special case of a generic operator $O_P$ specified in terms of a predicate $P(y)$ as

$$O_P A = \{y \in A \mid P(y)\}. \tag{3.1}$$

We recognize this as the Z-F Axiom of Comprehension. So $O_P A$ is a bona fide set. It follows that

$$x \subseteq y \Rightarrow O_P x = x \cap O_P y. \tag{3.2}$$

This relation holds in particular for $O_P$ taken equal to $\mathcal{R}$.

The Russell paradox is no longer relevant as a paradox. It is replaced by the operator $\mathcal{R}$ as examination of the proof of the following theorem reveals (compare (2.1)).

**Theorem 3.2:** $\forall A$, $\mathcal{R} A \notin A$. Moreover, $\mathcal{R} A$ is normal.

**Proof:** Assume there exists a set $z$ such that $\mathcal{R} z \in z$. Then by the definition of $\mathcal{R}$ there are two options, both of which lead to contradictions. Namely,

1. $\mathcal{R} z \in \mathcal{R} z$, in which case $\mathcal{R} z \notin \mathcal{R} z$,
2. $\mathcal{R} z \notin \mathcal{R} z$, in which case $\mathcal{R} z \in \mathcal{R} z$.

A corresponding result is

**Proposition 3.3:** $\forall A$, $A \notin \mathcal{R} A$.

**Proof:** By definition, if $x \in \mathcal{R}A$, then it is both true that $x \notin x$ and $x \in A$. Then $A \in \mathcal{R}A$ implies both $A \notin A$ and $A \in A$, a contradiction. □

We make the following observations associated with Theorem 3.2.

a) The collection of all sets is not itself a set.
b) Every set has an inside and an outside,
    where the inside of a set consists of its elements.
c) The complement of a set (the class of sets not in the given set) is not a set.
d) If $\forall y \in B, y \notin y$, then $B \notin B$.
e) If $\forall y \in C, y \in y$, we cannot conclude that $C \in C$.

S*ince $\mathcal{R}$ takes a part of A outside itself, note the relevance of b) to the ability to observe a set from the outside, a feature described in Sect. 1. To illustrate e) we first introduce the notion of the dual of a set.*

**Definition 3.4:** The dual $x^*$ of the set $x$ is given by

$$x^* = \{x^*, x\} \tag{3.3}$$

Existence and uniqueness of the dual of a set follows from the AFA. e) is illustrated by the following two examples.

Example 1: Since $\Omega \in \Omega$, taking $C = \Omega$ satisfies the hypotheses, and we have $C \in C$.
Example 2: Take $x$ and $y$ to be unequal normal sets, and let $C = \{x^*, y^*\}$. Then it Is easy to see that $C \notin C$.

Let $\mathcal{U} = \{x \mid x = x\}$ be the *class of sets*. ($\mathcal{U}$ is also referred to as the *universe of sets*.) Let $\mathcal{N} = \{x \mid x \notin x\}$ be the class of normal sets, and let $\mathcal{A} = \{x \mid x \in x\}$ be the class of abnormal sets. Then we have the following proposition concerning the classes $\mathcal{A}$, $\mathcal{N}$ and $\mathcal{U}$ and the Russell operator $\mathcal{R}$.

**Proposition 3.5:** a) $\mathcal{N}$ is a proper class.
   b) $\mathcal{U}$ is a proper class.
   c) $\mathcal{A}$ is a proper class.
   d) $\mathcal{U} = \mathcal{N} \cup \mathcal{A}$.
   e) $\mathcal{R}A = \mathcal{N} \cap A, \forall A$.

**Proof:** We shall prove a), b) and c). (See footnote 3).

a) Assume not. Then $\mathcal{N} = A$ for some set $A$. $\mathcal{R}A \notin A$ by Theorem 3.2. But then $\mathcal{R}A \in \mathcal{N}$, a contradiction.

b) Assume not. Then $\mathcal{U} = B$ for some set $B$. Now $\{x \in B \mid x \notin x\}$ is a set by the Axiom of Comprehension. However $\{x \in B \mid x \notin x\} = \mathcal{N}$ by definition. This is a contradiction since $\mathcal{N}$ is a proper class.

c) Suppose to the contrary that $\mathcal{A}$ is a set. Then there exists a unique a set $A$ such that $x \in A \Leftrightarrow x \in x$. Then using AFA, we see that $\forall y \in \mathcal{N},\ y^* \in A$. Let $C = \{x \in A \mid \exists y \in \mathcal{N} \text{ such that } x = y^*\}$. $C$ is itself a set (Axiom of Comprehension) that we can also write as

$$C = \{y^* \mid y \in \mathcal{N}\ \} = \{\{y^*, y\} \mid y \in \mathcal{N}\ \}. \tag{3.4}$$

Then using the Z-F Axiom of Union, we can write

$$\bigcup C = \bigcup_{y \notin y} \{y^*, y\}, \tag{3.5}$$

where $\bigcup$ on the left is the monadic union operator[4]. Then

$$\mathcal{R}(\bigcup C) = \bigcup_{y \notin y} \{y\} = \mathcal{N}. \tag{3.6}$$

This is a contradiction, since $\mathcal{R}(\bigcup C)$ is a set and $\mathcal{N}$ is a proper class. $\square$

### 3.2 Syntax

#### 3.2.1 Fundamental operators

We shall employ the following four basic dyadic set operations $\circ, \cup, \cap, -$, defined as follows.

$$\begin{aligned}
\circ: &\ (O_1 \circ O_2) x = O_1 O_2 x \\
\cup: &\ (O_1 \cup O_2) x = (O_1 x) \cup (O_2 x) \\
\cap: &\ (O_1 \cap O_2) x = (O_1 x) \cap (O_2 x) \\
-: &\ (O_1 - O_2) x = (O_1 x) - (O_2 x)
\end{aligned} \tag{3.7}$$

The last, the difference of operators, is defined in terms of set subtraction, given by the following Boolean rule.

$$x - y = x - (x \cap y). \tag{3.8}$$

The associative law $(O_1 O_2) O_3 = O_1 (O_2 O_3)$ follows from the definition of $\circ$.

To supplement $\mathcal{R}$ we introduce four additional basic operators $\mathcal{I}, \mathcal{E}, \mathcal{B}$ and $\mathcal{D}$, where

---

[4] The monadic union operator $\bigcup$ is defined as follows: $\bigcup A = \{x \mid x \in a \text{ for some } a \in A\}$, which is a set by virtue of the Axiom of Union.

    a) $I$ is the *identity operator*, $Ix = x$,
    b) $E$ is the *elimination operator*, $Ex = \emptyset$,
    c) $B$ is the *singleton operator*, $Bx = \{x\}$, and
    d) $D$ is the *duality operator*, $x^* = Dx = \{x^*, x\}$. (See Def. 3.4.)

While we defer introduction of semantics for the basic operators until Sect. 5, we make the following observations about them.

  *i)* $IO = OI$, for any operator O.
  *ii)* $E$ is not a right-zero operator, since for example, $BE \neq E$. Note that $E$ is *idempotent* ($E^2 = E$). Note also that $(BE)x = B\emptyset$, so that in particular, $(B^n E)x = B^n \emptyset$ for any non-negative integer $n$.
  *iii)* Since $\{x, y\} = \{x\} \cup \{y\}$, showing that $\{x, y\}$ is a derivative notion (see a) in Sect. 2.2), we can write

$$x^* = (Bx^*) \cup Bx. \tag{3.9}$$

The operators $B$ and $D$ are related as follows.

$$D = (BD) \cup B. \tag{3.10}$$

### 3.2.2 Properties of the basic operators

The quadruple of basic operators $B, E, I$ and $R$ and form a non-closed system illustrated in the following operator multiplication table.

|   | $E$ | $I$ | $R$ | $B$ |
|---|---|---|---|---|
| $E$ | $E$ | $E$ | $E$ | $E$ |
| $I$ | $E$ | $I$ | $R$ | $B$ |
| $R$ | $E$ | $R$ | $R$ | $RB$ |
| $B$ | $BE$ | $B$ | $BR$ | $BB$ |

**Table 3.1:** Multiplication table for the basic operators.

Next we introduce the *counter-Russell* operator, $T = I - R$. Note that

$$TA = A - RA = A \cap A. \tag{3.11}$$

Consider the following proposition relating $R$ and $T$ to normal and abnormal sets.

**Proposition 3.6:** Let $B$ be a normal set and $C$ an abnormal set. Then $RBB = BB$, but $RBC = \emptyset$. Alternatively, $TBC = BC$ but $TBB = \emptyset$.

We also have the following proposition exhibiting properties of $R$ and $B$.

**Proposition 3.7:** *i)* $x \in A \Leftrightarrow \mathcal{B}x \subset A$.
   *ii)* $\mathcal{R}A \notin A \Leftrightarrow \mathcal{B}\mathcal{R}A \not\subset A$.
   *iii)* $\mathcal{R}A \notin \mathcal{R}A$.

By Def. 2.2, $\mathcal{R}A$ is normal. Additional syntactic relations (conceptual operator statements) are given in the following theorem, proposition and corollary.

**Theorem 3.8:** a) $\mathcal{I} \cap \mathcal{R} = \mathcal{R}$.
   b) $\mathcal{B} \cap \mathcal{R} = \mathcal{E}$.
   c) $\mathcal{I} \cap (\mathcal{B}\mathcal{R}) = \mathcal{E}$.
   d) $\mathcal{I} \cap (\mathcal{R}\mathcal{B}) = \mathcal{E}$.
   e) $\mathcal{R}\mathcal{B} = \mathcal{B} - \mathcal{I}$.
   f) $\mathcal{R}\mathcal{B} - \mathcal{B}\mathcal{R} \neq \mathcal{E}$.
      $\mathcal{B}\mathcal{R} - \mathcal{R}\mathcal{B} \neq \mathcal{E}$.

The two statements in f) are not the same since there is no monadic minus for sets. Examples illustrating the two relations in f) are $(\mathcal{R}\mathcal{B} - \mathcal{B}\mathcal{R})\mathcal{B}\mathcal{D}\emptyset \neq \emptyset$ and $(\mathcal{B}\mathcal{R} - \mathcal{R}\mathcal{B})\Omega \neq \emptyset$, respectively.

**Proposition 3.9:** $\mathcal{R}\mathcal{B}\mathcal{R} = \mathcal{B}\mathcal{R}$.

**Proof:** The proof is a direct consequence of Prop. 3.6 and Prop. 3.7. □

**Corollary 3.10:** $(\mathcal{R}\mathcal{B} - \mathcal{B}\mathcal{R})\mathcal{R} = \mathcal{E}$, and $(\mathcal{B}\mathcal{R} - \mathcal{R}\mathcal{B})\mathcal{R} = \mathcal{E}$.

This corollary gives a connection between the Prop. 3.9 and relation f) in Theorem. 3.8.

### 3.2.3 Characterization of $\mathcal{R}$

The following proposition and corollary gives a complete characterization of $\mathcal{R}$.

**Proposition 3.11:** If
$$x \subseteq y \Rightarrow Ox = x \cap Oy, \tag{3.12}$$
then $O\mathcal{B}$ uniquely determines $O$.

The conclusion of the proposition may be restated alternatively as

$$\forall x, \, Ox = \{y \in x \mid O\mathcal{B}y = \mathcal{B}y\}. \tag{3.13}$$

**Proof of Prop. 3.11:** We make the following two preliminary observations. (i) The hypothesis implies that $\forall x, \, Ox \subseteq x$, and hence that (ii) $O^2x = Ox \cap Ox = Ox$. Continuing we now address the question: when is $y \in Ox$? However $y \in Ox \Leftrightarrow \mathcal{B}y \subseteq Ox$, by definition. In the hypothesis we may replace $x$ with $\mathcal{B}y$ and $y$ with $Ox$ to conclude that

$$O\mathcal{B}y = \mathcal{B}y \cap O(Ox) = \mathcal{B}y \cap Ox, \tag{3.14}$$

the last employing (ii). Now $By \subseteq Ox \Leftrightarrow By \cap Ox = By$ by definition. This and (3.13) implies $y \in Ox$ if and only if

$$OBy = By. \tag{3.15}$$

From this and (i) we conclude that $y \in Ox \Leftrightarrow y \in x$ and $OBy = By$. □

**Corollary 3.12:** Let the operators $O_1$ and $O_2$ satisfy the hypothesis of the proposition, and let $O_1 B = O_2 B$. Then $O_1 = O_2$.

The notion of an operator called a selector (compare Def. 3.1 *f*) is specified as follows.

**Definition 3.13:** A class of operators called *selectors*[5] are those that satisfy the hypothesis $x \subseteq y \Rightarrow Ox = x \cap Oy$ of Prop. 3.11.

Selectors form a commutative system, as the following proposition shows.

**Proposition 3.14:** If $O_1$ and $O_2$ are selectors, then $O_2 O_1 = O_1 O_2 = O_1 \cap O_2$.

**Proof:** Let $z = O_1 x \subseteq x$. Then $O_2 z = z \cap O_2 x$. Then $O_2(O_1 x) = (O_1 x) \cap O_2 x = O_1 x \cap O_2 x = (O_1 \cap O_2) x$. □

In particular, consciousness operators $\mathcal{K}$ to be introduced in Sect. 3.3, being selectors, commute. Theorem 3.8(f) provides an example of a non-commuting pair of operators.

Characterization of the Russell operator is the subject of the following theorem.

**Theorem 3.15:** $\mathcal{R}$ is characterized by the following two properties.

1. $x \subseteq y \Rightarrow \mathcal{R}x = x \cap \mathcal{R}y$. (3.16)
2. $\mathcal{R}B = B - I$.

**Proof:** #1 follows from (3.2). #2 is the result e) of Theorem 3.3. Then the proof of the characterization $\mathcal{R}$ of is an immediate consequence of Prop. 3.11 and Corollary 3.12.

### 3.2.4 Schematic illustrating syntax of sets and operators

The Venn type diagram in Fig. 3.2 illustrates some of the notions being discussed. The diagram is intended to be composed in a homeomorphic representation of the Euclidean plane. In the diagram sets and classes are represented by open rectangles. That is, they do not contain their boundaries. For example in terms of the rectangular coordinates $\alpha$ and $\beta$ in the plane, the empty set is given by $\emptyset = \{\alpha, \beta \mid \alpha^2 + \beta^2 < 0^2\}$. The class of abnormal sets is represented by the largest shaded rectangle. The class of normal sets is represented by the largest un-shaded rectangle. A set's name is displayed at the tail of an arrow pointing to that set. To interpret the diagram, consider the large rectangle *A* in the middle of the figure. *A* is positioned in a general position, that is, so that some of its elements are abnormal (shaded) and some are normal. *F* is a subset of *A*. *C* (toward the lower right) is a set all of whose elements are normal, and *G* is a subset of *C*. *D* (upper right) is a set all of whose elements are abnormal, and *H* is a subset of *D*. The remaining sets indicate the result of applying one or more operators to the sets and subsets just identified.

---

[5] Thanks to the referee for noting that a notion of selector occurs in relational databases.

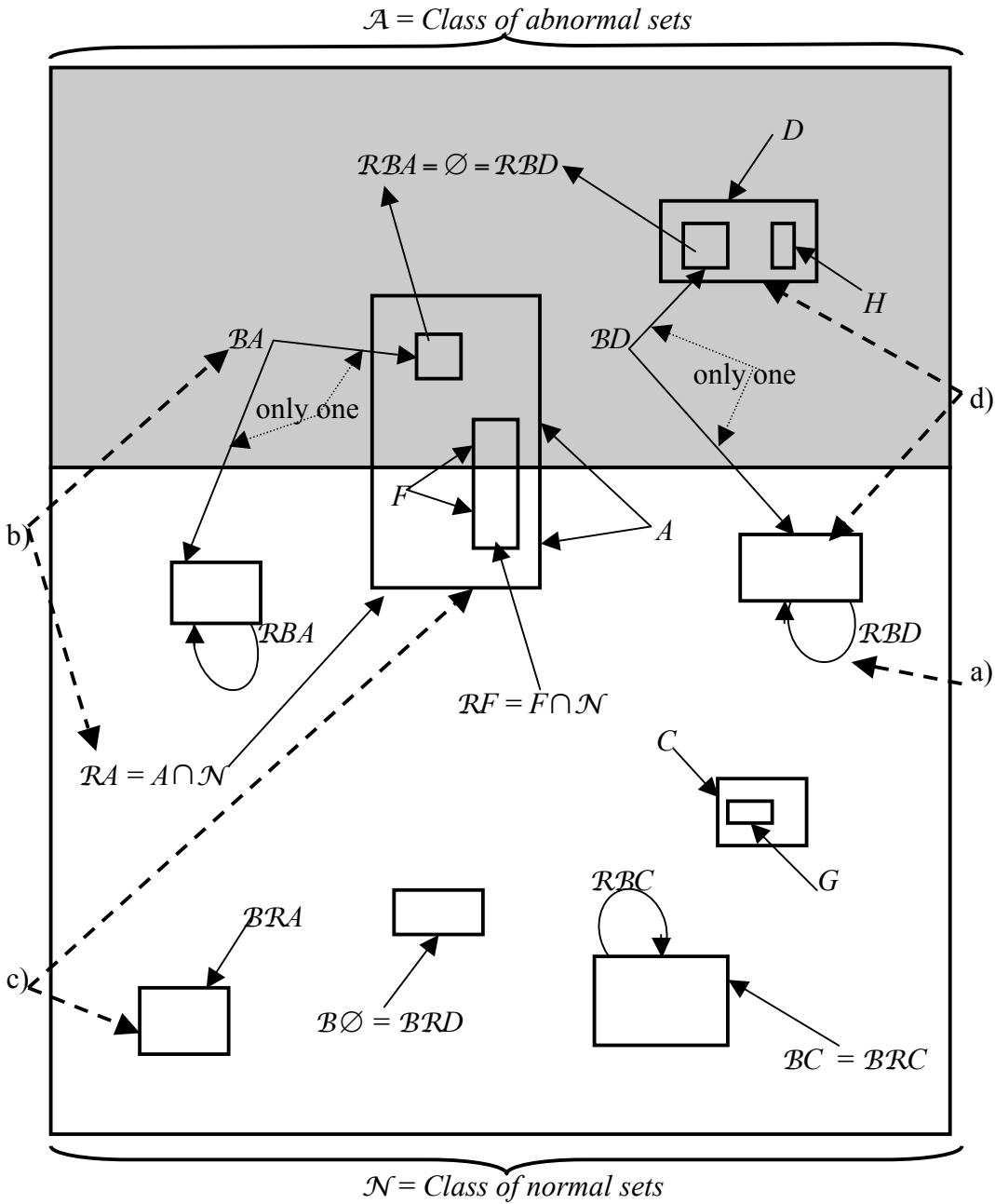

**Figure 3.2:** A schematic illustrating the properties assembled in Theorem 3.3

Illustrated in Fig. 3.2 are 6 possibilities for sets and 7 for fundamental operators:

2 for *A*, depending on whether $\mathcal{B}A \subseteq A$ or not. (See the phrase "only one" in the figure.)
1 for set *C*, namely, $\mathcal{B}C \not\subseteq C$.
2 for *D*, depending on whether $\mathcal{B}D \subseteq D$ or not. (See the phrase "only one" in the figure.)
1 for $\varnothing$, a technical possibility, since $\varnothing$ can not be illustrated.
The 7 illustrated fundamental operators are $\mathcal{E}, \mathcal{I}, \mathcal{B}, \mathcal{R}, \mathcal{T}, \mathcal{BR},$ and $\mathcal{RB}$, although $\mathcal{E}, \mathcal{I}$ and $\mathcal{T}$ are illustrated implicitly.

The conclusions a) - d) of Theorem 3.8 are illustrated in the figure by the sets and/or labels of sets that are pointed to by dashed arrows with the corresponding labels. These labels are placed in the margins of the figure. For example, the c) in the left hand margin labels both a dashed arrow pointing to the set $\mathcal{BR}A$ and a dashed arrow pointing to the label of the set *A*. These two sets, shown as disjoint in the figure, illustrate conclusion c) of the theorem. One can see that conclusions e) and f) are also illustrated. The result

$$x \subseteq y \Rightarrow \mathcal{R}x = x \cap \mathcal{R}y, \qquad (3.17)$$

following from (3.2) showing that $\mathcal{R}$ is a selector is illustrated in its three different cases.

1. $\mathcal{R}F$, the part of *F* in $\mathcal{N}$ equals $F \cap \mathcal{R}A$.
2. $H \subseteq D \in \mathcal{A}$ then $\mathcal{R}H = \varnothing$.
3. $G \subseteq C \in \mathcal{N}$ then $\mathcal{R}G = G \cap C = G$.

### 3.3 Semantics and consciousness operators

We now develop a model in which *experience* and *consciousness* are taken as primitives. These primitives may be composed of layers. In this case, our primitives model the corresponding basic layers, namely what we have knowledge and understanding about through our sensations and perceptions (this last being a Cantor-like statement). When necessary for clarity, the basic layers shall be called *primary experience* and *primary consciousness*, respectively. While we perceive these basic layers, they are essentially ineffable. The higher layers, should they exist, might very well be beyond ineffability. We focus on the basic layers, and we take our primitives to be models of them. Our goal to specify an illuminating axiomatic system for these primitives. So we may say that as with set theory, we commence with a Cantor-like (naïve) manner and then refine it by means of an axiomatic approach.

We shall characterize a collection of operators called *consciousness operators*, the generic element of which is denoted by $\mathcal{K}$. We take a set *x* to model a primary experience. Such a set, being a primitive, may be viewed as a Platonic object. Then our *Semantic Thesis* is stated as follows.

**Definition 3.16** (**Semantic Thesis**): *Consciousness is a result of a consciousness operator being applied to experience.*

We now give the first four axioms of an open (and developing) system that serves to characterize the experience and consciousness primitives. The axioms and their semantic interpretations justify the Semantic Thesis. We begin with the following definition.

**Definition 3.17:** Let $x$ model a primary experience. Then $\mathcal{K}x$ models the *awareness*, an *induced experience*. *Consciousness* is an instance of a specific operator $\mathcal{K}$ acting on experience.

The first three axioms along with their semantic interpretations and a name for each are displayed in the following table.

|    | Axiom | Semantic interpretation of the axiom | Name of Axiom |
|----|-------|--------------------------------------|---------------|
| a) | $\forall x,\ \mathcal{K}x \subseteq x$ | Experience generates its own awareness | Generation |
| b) | $\forall x,\ x \notin \mathcal{K}x$ | Awareness does not generate the primary experience | Irreversibility |
| c) | $\forall x,\ \mathcal{K}x \nsubseteq x$ | Awareness is removed from experience | Removal |

**Table 3.3:** First three axioms for a consciousness operator.

Axioms a) and b) are motivated by the properties of the Russell operator a) and b), respectively given in Theorem 3.8. By taking $y = \mathcal{K}x$ in the following Prop. 3.18 and using axioms a) and b), we conclude that a set of the form $\mathcal{K}x$ is normal. This may be interpreted semantically as the normality of awareness.

The following analytic statement of axiom c)

$$(\mathcal{BK}x) \cap x = \varnothing. \tag{3.18}$$

follows by noting that $y \notin x$ and $y \subseteq x$ implies $y \notin y$.

The following table displays algebraic statements of these axioms along with examples of operators that violate each statement. The existence of $\Omega$ shows that $\mathcal{B}$ and $\mathcal{I}$ violate c).

|    | Algebraic statement | Violating examples |
|----|---------------------|--------------------|
| a) | $O \cap I = O$ | $O = \mathcal{B}$ |
| b) | $\mathcal{B} \cap O = \mathcal{E}$ | $O = \mathcal{B}, \mathcal{T}$ |
| c) | $(\mathcal{B}O) \cap I = \mathcal{E}$ | $O = \mathcal{E}, \mathcal{B}, \mathcal{I}$ |

**Table 3.4:** Algebraic statements of the first three axioms for a consciousness operator.

We now append a fourth axiom that in fact is stronger than axiom a).

| d) | If $x \subseteq y$, then $\mathcal{K}x = x \cap \mathcal{K}y$ | Awareness of a sub-experience is determined by the sub-experience and awareness of the primary experience | Selection |
|---|---|---|---|

**Table 3.5:** The fourth consciousness operator axiom.

Axiom d) is motivated by its syntactic counterpart expressed by the condition (3.12) of Prop. 3.11. Axiom d) is the statement that $\mathcal{K}$ is a selector.

The consistency of the axioms a) - d) is demonstrated by producing an operator that satisfies all of them. Indeed, $\mathcal{R}$ is such an operator as the following theorem shows.

**Theorem 3.18:** The Russell operator $\mathcal{R}$ satisfies the axioms a), b), c) and d).

**Proof:** The proof follows from properties of $\mathcal{R}$ assembled in Sect. 3.1 and Sect. 3.2.

The following result describes the action of $\mathcal{K}$ on the primary experience $\mathcal{B}\emptyset$.

**Proposition 3.19:** $\mathcal{KBE} = \mathcal{BE}$.

**Proof:** Axiom c) implies that $\mathcal{KBE}\emptyset \neq \mathcal{BE}\emptyset$. Axiom a) implies that $\mathcal{KB}\emptyset \subseteq \mathcal{B}\emptyset$. Together these two statements imply that $\mathcal{KB}\emptyset = \mathcal{B}\emptyset$. □

There are other operators besides $\mathcal{R}$ that satisfy axioms a) - d), as the following example $C$ of a consciousness operator shows.

$$Cx = \{y \in x \mid y \notin y; \forall z \in y, z \neq \Omega\}. \tag{3.19}$$

So since all elements of $Cx$ are normal, $Cx$ is a normal set. Moreover $Cx \subseteq \mathcal{R}x$, so that $C$ is a sub-operator of $\mathcal{R}$. To show that $C \neq \mathcal{R}$ note that for the set $A = \{\{\emptyset, \Omega\}\}$, we have $\mathcal{R}A = A$, but $CA = \emptyset$. To show that $C$ satisfies the axioms, we proceed as follows.

a) By definition $Cx \subseteq x$, so axiom a) is satisfied.
b) Since $Cx \subseteq \mathcal{R}x$ and $x \notin \mathcal{R}x$, then $x \notin Cx$. So axiom b) is satisfied.
c) To prove that $C$ satisfies axiom c), we show the algebraic equivalent axiom c). Namely that $\mathcal{BC} \cap \mathcal{I} = \mathcal{E}$. Then suppose $\exists z$ such that $Cz \in z$. There are two options.
   1. $Cz \in Cz$. This implies that $Cz \notin Cz$, a contradiction, since by definition every element of $Cx$ is normal.
   2. $Cz \notin Cz$. This implies that either $Cz \in Cz$ or $\Omega \in Cz$. Hence $Cz \notin Cz$ implies $\Omega \in Cz$. However $\Omega \in \Omega$, contradicting the normality of $Cz$.
d) (3.2) shows that $C$ satisfies axiom d).

While the axioms of a consciousness operator appear to be limiting, we are able to exhibit an infinite collection "$\{\mathcal{K}_A | A \in \mathcal{U}\}$" of such operators. In particular generalizing (3.19) yields the following.

$$\mathcal{K}_A x = \{y \in x \,|\, y \notin y,\, \mathcal{T}(y \cap A) = \varnothing\}. \tag{3.20}$$

Comments on a connection of qualia to a diagonalization of $\mathcal{K}_A$ are given in Sect. 6.

## 4. Labeling of Graphs, Histogram Construction, M-Z Equation, Neural Networks

We begin with a prescription for labeling a collection (Sect. 2.1). This is extended to a technique for decorating a labeled graph. Given a graph, this procedure forms the basis for inducing the existence of a set intrinsically associated with a labeled graph. Then a construction of what we call a histogram is made. The latter is a novel tool used in proposing the M-Z equation, which comes from a synthesis of Aczel's theory of decorating labeled graphs and the theory of neural networks. An interpretation is made that portrays the sets decorating a labeled graph as Platonic constructs. So this application and interpretation constitute a theory of consciousness constructed on the foundations developed here.

### 4.1 Labeling of graphs
A *labeling* $\lambda$ of $\Gamma$ is a set valued function of the nodes $N$ of $\Gamma$.

$$a \mapsto \lambda a,\, \forall a \in N. \tag{4.1}$$

A decoration of a labeled graph is a set valued function $a \mapsto d_\lambda a$, where (compare (2.2))

$$d_\lambda a = \{d_\lambda b \,|\, a \to b\} \cup \lambda a,\; \forall a \in N. \tag{4.2}$$

This system of equations along with the following theorem shows how labeled decorations are a basis for inducing the existence of a set intrinsically associated with the labeled graph. (Compare the notion of the picture of a graph in Sect. 2.2.) Existence and uniqueness of $d_\lambda$ is the subject of the following theorem.

**Theorem 4.1:** Given $\Gamma = (N, E, \lambda)$, a corresponding decoration $d_\lambda$ of $\Gamma$ exists and is unique. (Aczel, Theorem 1.10.) (Compare with AFA in Sect. 2.2.)

**Example:** Take $\Gamma = (N, E)$ to be the graph of the set $\Omega$. $\Gamma$ is specified by $N = \{a\}$ and $E = \{a \to a\}$. Then with $\lambda a$ being any set, we have

$$d_\lambda a = \{d_\lambda a\} \cup \lambda a. \tag{4.3}$$

If $\lambda a = \{b\}$, a singleton, then $d_\lambda a = \{d_\lambda a, b\}$. Then $d_\lambda a = b^* = \mathcal{D}b$ is the dual of $b$.

## 4.2 The histogram construction

We now introduce a construct called the histogram of a function that replaces a set valued function on a collection by a set valued function on a pure set. The construct is used to apply Theorem 4.1 to a collection of graphs abstracted from models of brain circuitry to be introduced in Sect. 4.4.

Let $A$ be a collection of unknown elements, and let $B$ be a set. Consider a mapping, $f : A \to B$, and define

$$f^{-1}(b) = \{a \in A \mid f(a) = b\}, \forall b \in B. \tag{4.4}$$

We suppose that the number of elements in this set, $|f^{-1}(b)|$ is finite for every $b$. Then the histogram of a mapping is specified as follows.

**Definition 4.2 (Histogram):** The *histogram* $H_f$ of $f$ is the following set of ordered pairs.

$$H_f = \{(b, |f^{-1}(b)|) \mid b \in B, f^{-1}(b) \neq \emptyset\}. \tag{4.5}$$

Note that $H_f$ is a bona fide set (see Sect 2.1), and in particular, that $H_f \subseteq B \times \mathbf{N}_+$.

## 4.3 The M-Z equation, the weight function, the voltage function

We call a function $w : E \to \mathbf{Q}$, a *weight function*. The rationals $\mathbf{Q}$ comprise a set, since each non zero rational number $q$ corresponds to the triple $(m, n, \pm)$, where $\pm m/n$ ($m$ and $n$ being relatively prime natural numbers) is the value of $q$. The choice of the rationals for the range of $w$ is made for definiteness and clarity.

Let $E_a$ denote the set of edges of a graph $\Gamma$ that terminate in the node $a$, that is,

$$E_a = \{(p, a) \mid p \in N, p \to a\}, \forall a \in N. \tag{4.6}$$

We make the following local finiteness hypothesis.

**Hypothesis 4.3:** $\forall a, E_a$ is finite.

Then let[6] $w_a = w|_{E_a}$, so that $w_a : E_a \to \mathbf{Q}$ is a function from a finite unordered collection into the rationals.

Let $H_{w_a}$ be the histogram of the mapping $w_a$. $H_{w_a}$ is a finite set since $E_a$ is a finite collection. Note that

$$H_{w_a} \subseteq \mathbf{Q} \times \mathbf{N}_+. \tag{4.7}$$

We now define the M-Z equation.

---

[6] Recall that $f|_y$ denotes the restriction of a mapping $f$ to a sub collection $y$ of the domain of $f$.

**Definition 4.4 (M-Z equation):** Given $\Gamma = (N, E, w)$, label $\Gamma$ with the labeling $\lambda : a \mapsto H_{w_a}$. Then the labeled decoration of $\Gamma$ is specified by the M-Z *equation*, namely

$$d_\lambda a = \{d_\lambda b \,|\, a \rightarrow b\} \cup H_{w_a}, \quad \forall a \in N. \tag{4.8}$$

Comparing (4.8) to (2.2) where a decoration is defined, we may interpret the set $H_{w_a}$ as a forcing term in the M-Z equation for the decoration $d_\lambda a$. In a forthcoming work (Miranker, Zuckerman 2008), a number of examples and applications of the M-Z equation is assembled.

We shall be interested in an extension of the above development that involves what we call a *voltage function* $v : N \rightarrow \{0,1\}$. (The choice of $\{0,1\}$ is made for definiteness and clarity.) Take

$$E_{a,v} = \{(p, a) \,|\, p \rightarrow a, v(p) = 1\}, \forall a \in N, \tag{4.9}$$

and let $w_{a,v} = w|_{E_{a,v}}$. Note that the histogram $H_{w_{a,v}} = \varnothing$ if $E_{a,v} = \varnothing$. Now label $\Gamma$ with the labeling $\lambda : a \mapsto H_{w_{a,v}}$. Then the M-Z equation that specifies the labeled decoration of $\Gamma = (N, E, w, v)$ is given by (4.8) with $H_{w_a}$ replaced by $H_{w_{a,v}}$. Namely,

$$d_\lambda a = \{d_\lambda b \,|\, a \rightarrow b\} \cup H_{w_{a,v}}, \quad \forall a \in N. \tag{4.10}$$

In Fig. 4.1 we give a schematic of the process of labeling a neuron with a histogram

*r active efferent neurons* (*each with* $v = 1$)

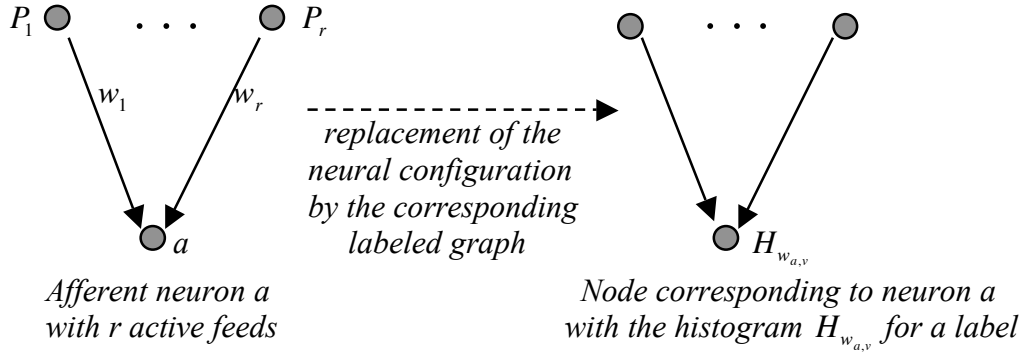

*Afferent neuron a*              *Node corresponding to neuron a*
*with r active feeds*              *with the histogram* $H_{w_{a,v}}$ *for a label*

**Figure 4.1:** Labeling a neuron with a histogram

Recall that the histogram is constructed ignoring the numbering of the efferents.

## 4.4 Application to a neural network model of brain circuitry

The brain is commonly taken as the seat of consciousness, the latter supervening on the workings of the brain's neural networks. (While for some, it is the entire physical body and even the environment that is taken as the seat of consciousness, there is no loss of meaning for our argument to take the more limited view.) We shall show how our constructs apply to a neural network to produce a labeled decorated graph. This in turn allows us to incorporate neural networks into the mathematical foundation of consciousness. Take a neuron $Q$ and trace its inputs (afferents) backward and its outputs (efferents) forward to elaborate a neural network. Replacing a neuron and its dendritic and axonal processes by a node and its synapses by directed edges, there results a graph $\Gamma$ emanating from the node (also called $Q$) corresponding to the chosen neuron. Typically this network has reentrant connections, and so, $\Gamma$ is non-well-founded. An illustration of a possible $\Gamma$ fragment is given in Fig. 4.2. Note the correspondence to the cords and knots of Kanger, 1957.

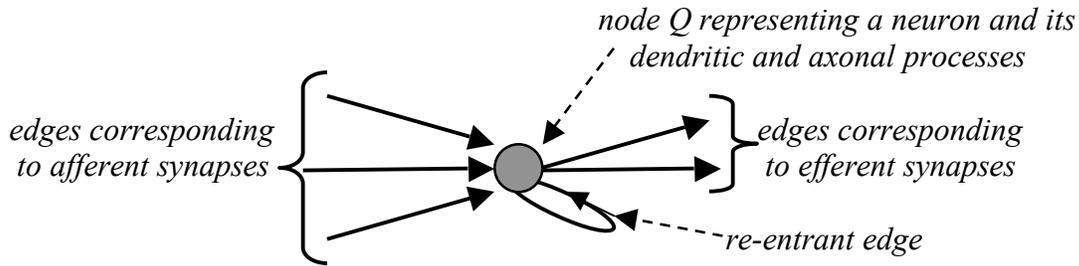

**Figure 4.2:** A neural net with a single node $Q$ interpreted as a graph.

This network and so also $\Gamma$ is associated with two families of parameters, namely, the synaptic weights $w$ and the output of its neurons' activities. The latter are expressed as voltages, denoted $v$. Hebb's rule is the customary model of synaptic weight change. The changes in voltage outputs are modeled by input-output threshold equations, the simplest version of which is the McCulloch-Pitts model (Haykin 1999). For clarity, we employ the simplest meaningful form of these two models, using them to specify updates of $w$ and $v$, the latter written as $w_{old}(a \to b) \xrightarrow{update} w_{new}(a \to b)$ and $v_{old}(a) \xrightarrow{update} v_{new}(a)$, respectively.

**Definition 4.5 (Hebb's rule):**

$$w_{new}(a \to b) - w_{old}(a \to b) = \alpha v_{old}(a) v_{new}(b). \qquad (4.11)$$

In (4.11), $a$ is an efferent neuron, $b$ is a corresponding afferent and $w(a \to b)$ is the weight of the synapse connecting neuron $a$ to neuron $b$. (For convenience we allow at most one such connection per pair of neurons.) The voltage $v(a)$ is the neuronal activity of $a$. For consistency with Sect. 4.3, the scaling constant $\alpha$ is chosen to be a rational number.

**Definition 4.6 (McCulloch-Pitts equation):**

$$v_{new}(a) = h\left(\sum_{p:p\to a} w_{old}(p\to a)v_{old}(p) - \theta\right). \quad (4.12)$$

In (4.12) $h$ is the Heaviside function, the real number $\theta$ is a threshold, and the sum is over all neurons $p$ that forward connect directly to neuron $a$.

Note that (4.11) and (4.12) form a coupled dynamical system.

**4.4.1 Neural state, its decoration. The neural net semantic thesis**

At any instant of time, the coupled dynamical system (4.11) and (4.12) may be viewed as specifying the current states of the functions $v$ and $w$. We use the term *neural state* to describe this instantaneous state of the neural assembly. Referring to the weight and voltage functions of Sect.4.3, we use the $v$ and $w$ to specify a labeling, $\lambda_{w,v}: a \mapsto H_{w_{a,v}}$ of the graph $\Gamma$ as described in that section. Then we may use (4.10) to specify a labeled decoration, $d_{\lambda_{w,v}}$ of $\Gamma$. We shall also refer to $d_{\lambda_{w,v}}$ as the labeled decoration of the corresponding neural state. We now state our *Neural Net Semantic Thesis*. (See the Semantic Thesis of Def. 3.16.)

**Definition 4.7 (Neural Net Semantic Thesis):** *Each value of the Platonic function $d_{\lambda_{w,v}}$ encodes a dynamic preconscious experience associated with the corresponding neuron (i.e., node of $\Gamma$).*

As the brain processes information, the weights and voltages change as characterized by the Hebbian dynamics and the Mc-P dynamics. These in turn inform changes in associated preconscious experiences.

**4.4.2 Platonism**

Neural networks are physical, that is, they may be observed and their weights and voltages can be measured. The set values of the labeled decorations $d_{\lambda_{w,v}}$ are not physical. *Since they are located in some virtual space, we regard a value of $d_{\lambda_{w,v}}$ as Platonic.* (Compare Schrödinger's quote in Sect. 1.)

If $\Gamma$ is well-founded, its labeled decoration can be constructed in a recursive manner (Aczel 1988). However while the AFA supplies an existence statement for the decoration of a non-well-founded graph, it does not give a method to construct that decoration. The universe of graphs is divisible into two parts, one in which labeled decorations are recursively computable and the complement. The computability for graphs in the first part is a reason for classifying these corresponding sets as physical and not Platonic. The non-computability of graphs in the second part reinforces their Platonic status.

### 4.4.3 Example: A neural state instantiating a concept; Memes and themata

Consider the model neural network in Fig. 4.3a composed of three McCulloch-Pitts neurons, $a$, $b$, $c$ with the synaptic weights $w_{ba}$, $w_{ca}$, $w_{bc}$ (where for example, $w_{ba}$ denotes $w(a \rightarrow b)$) and with the voltages $v(a) = v(c) = 0$ and $v(b) = 1$. With these data and with the time frozen, the network becomes what we have called a *neural state*. When the APG (shown in Fig 4.3b) associated with this neural state is appropriately labeled with the specified weight and voltage data and then decorated, the diagram in Fig. 4.3b, a picture of a particular set $\Theta$ results. Since the source voltages vanish, the histograms are empty. Then for the sets of the decoration, we have $\Theta = \{B,C\}$, where $B = \varnothing$ and $C = \{\varnothing\}$.

Diagram Fig 4.3b, illustrating a decorated labeled APG, arises from the neural state in Fig. 4.3a. The APG in Fig. 4.3b is a representation of the Von Neumann ordinal 2, so that APG *is an instantiation of the ordinal 2*. (A number of additional examples are found in Miranker, Zuckerman 2008, where concepts and their instantiations are termed memes, and where the instantiation of an interpretation of a concept as a set is called the thema of that concept (of that meme).) We see that the set $\Theta$ decorating the point is the thema of the meme instantiated by this APG. The thema $\Theta$ along with the diagram in Fig. 4.3a are Platonic instantiations. The corresponding actual neural state being modeled (as by the model in Fig. 4.3a) is physical instantiation of that meme. The neural state (as illustrated in Fig. 4.3a) and the APGs (such as illustrated in Fig. 4.3b) are only examples of a vast number (in principle an unbounded number) of neural states and corresponding APGs that have this same thema $\Theta$. All such APGs are pictures of the set $\Theta$, so by analogy, we might say that each such meme whose thema is $\Theta$ is a picture of $\Theta$.

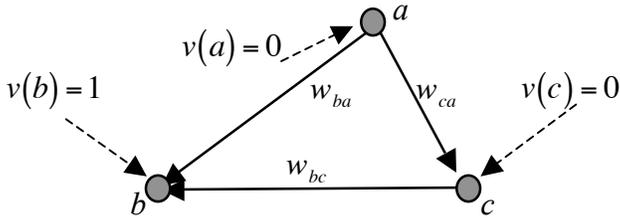
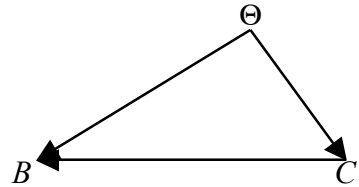

**Fig 4.3a:** Neural state with neurons $a$, $b$, $c$  
**Fig 4.3b:** Corresponding APG with point $\Theta$

### 4.5 Correspondence of the semantic theses, a neuro-physiological thesis

The relationship among the semantic theses of Sects. 3.3 and 4.4 is shown in Fig. 4.4. Each arrow in Fig. 4.4 describes a flow of information. The lowest arrow is a flow of physical information. The second is a flow from physical to Platonic information. The highest is a flow of psychic (Platonic) information. Fig 4.4 portrays the following notion.

**Definition 4.9 (Neuro-physiological thesis):** The *neuro-physiological thesis* (Fig. 4.4, lower left), denoted the movement of sensory information from a sense organ to the brain where it is processed to frame an internal physical representation of that information, and from where a primitive called consciousness is made manifest in a virtual space.

Note a parallel between the information flow in Fig. 4.4 with Plato's *line of knowledge*.

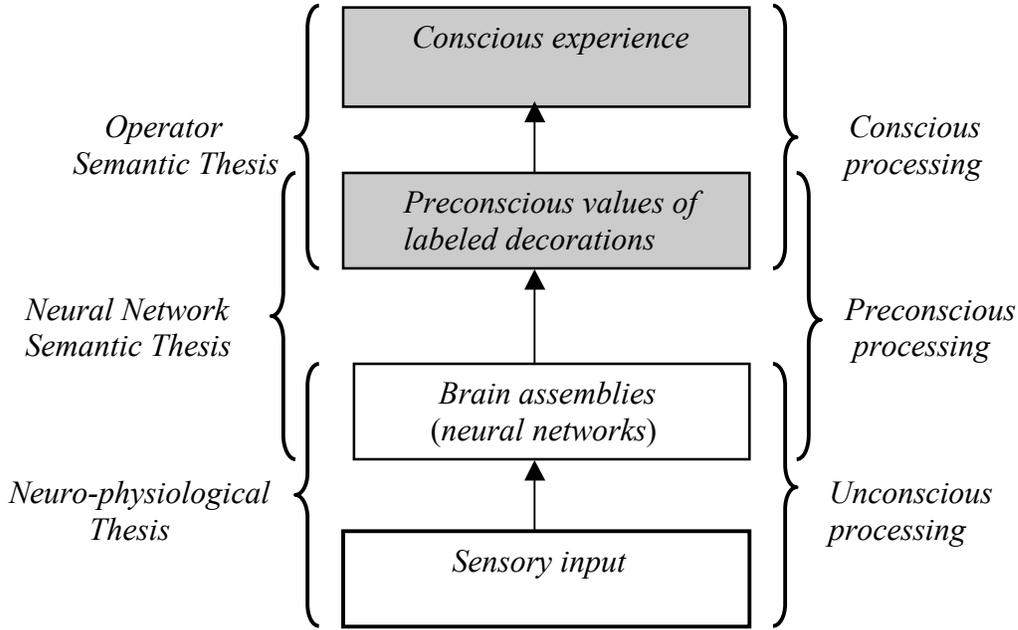

**Figure 4.4:** The theses of consciousness. The shading demarks the Platonic realm.

### 5. Observations: syntactic and semantic nomenclature

In this section we shall develop an elaboration of the consciousness operators introduced in Sect. 3.3 by exploiting aspects of the constructs developed in Sect. 4. The fundamental operators of set theory introduced in Sect. 3.2.1 will contribute as well. This elaboration will provide additional applications of the theory presented here. This is motivated by the syntactic and semantic nomenclature ascribed to these operators, which is summarized in the Table 5.1. Also shown in this table is an interpretation of each operator along with the axiom(s) that the operator codifies.

| Op | Syntactic | Semantic | Interpretation | Axiom(s) |
|---|---|---|---|---|
| $\mathcal{E}$ | Elimination | Erasing/Forgetting | Erases set representing Platonic experience | Existence of $\varnothing$ |
| $\mathcal{I}$ | Identity | Accepting/Receiving | Leaves set unchanged | Extension |
| $\mathcal{B}$ | Brace | Conceiving | Creates higher order set (a singleton) out of a set | Pair and singleton |
| $\mathcal{R}$ | Russell | Perceiving | Bifurcates set contents & retains normal elements | Comprehension |
| $\mathcal{T}$ | Anti-Russell | Rejecting/Denying | Counters $\mathcal{R}$, retaining the abnormal elements | Union |
| $\mathcal{D}$ | Duality | Reinforcing/Elaborating | Elaborates the concept of a set | AFA |

**Table 5.1:** Semantic interpretations of basic operators

In Fig. 5.2 we schematize the flow of information from sensory input to conscious experience. The upper boxes describe the syntactic level, the lower the semantic. A neural network in the brain typically corresponds to a non-well-founded graph. Hopfield networks supply examples. The corresponding labeled decorations are not recursively computable. They are schematized in the box labeled "*Functions $d_\lambda$ with values in virtual sets*" in Fig. 5.2. Is it a time dependent one of these deorations that emerges into consciousness? If so, how is the corresponding neural network selected?

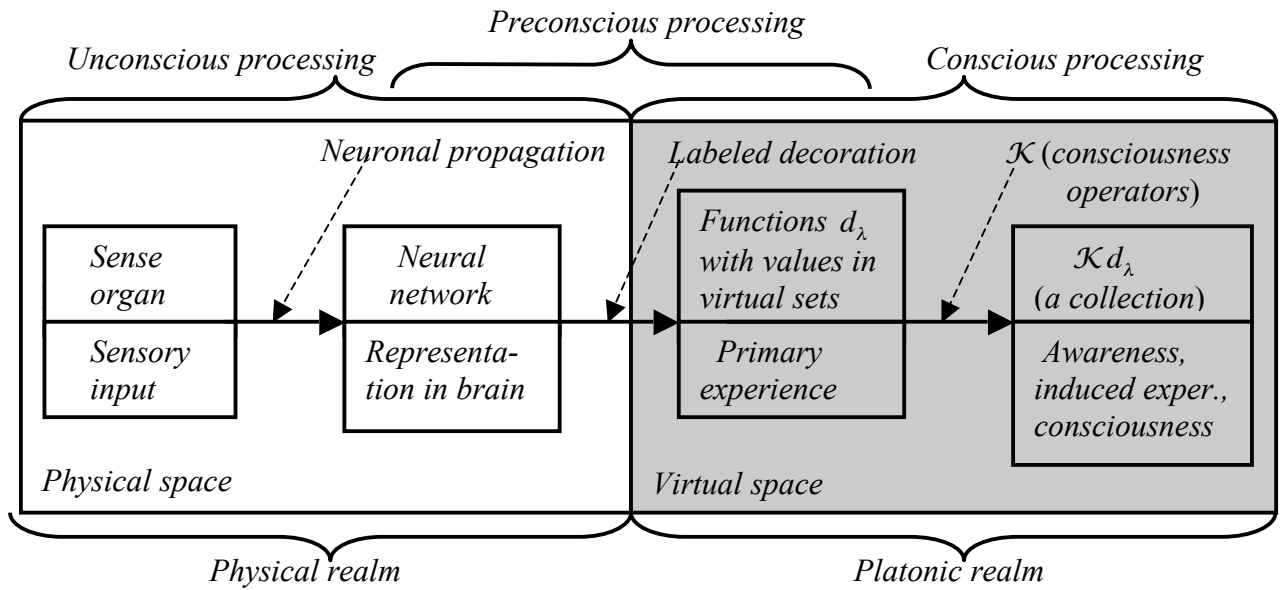

**Figure 5.2:** Consciousness: Syntactic and semantic views of the processing from the physical to the virtual. Shading distinguishes the ideal Platonic realm from the physical.

## 6. Future directions

**6.1 Applications of the M-Z equation**
In ongoing work (Miranker, Zuckerman 2008) numerous examples of applications of the M-Z equation are developed in the context of neural network modeling of brain circuitry. In particular the development of the notion of memes (i.e., concepts that are represented by pictures) and their themata (Platonic interpretations or themes of a concept collection) briefly introduced in Sect. 4.4.1 is made and complemented with examples. That work is extended to multi-graphs, a more faithful model of brain circuitry. Then a notion of histogram dynamics (Sect. 4.2) is introduced as a way to study the discrete time dependence of the associated notions of awareness and consciousness. Abstractions of those dynamics are developed as a tool for their study. The non-well-founded set theoretic framework that provides the context of these developments leads to a notion of a hierarchy of perceptual realities that we expect will inform understanding of the features of consciousness.

## 6.2 Diagonalization of the family of operators $\mathcal{K}_A$; qualia:

We propose to study the diagonalization $\mathcal{K}^{diag}$, of the family of operators $\mathcal{K}_A$ introduced in (3.20). The diagonalization is specified as follows.

$$\mathcal{K}^{diag} x = \{y \in x \mid y \notin y, \text{ and } \forall z \in x \cap y, z \notin z\}. \tag{6.1}$$

$\mathcal{K}^{diag}$ satisfies axioms a) – c) of Sect. 3.2. However taking $A_1 = \{\{\emptyset, \Omega\}\}$ and $A_2 = \{\emptyset, \Omega, \{\emptyset, \Omega\}\}$ as two choices for $A$ in (3.20), it follows that $\mathcal{K}^{diag} A$ is not a subset of $\mathcal{K}^{diag} B$. So failing axiom d) precludes $\mathcal{K}^{diag}$ from being a consciousness operator. Nevertheless we expect $\mathcal{K}^{diag}$ to be an operator of interest. For instance, take the set $d_\lambda a$ specified by the M-Z equation in (4.8) and put $a$ equal to $p$, the point of a graph, that graph corresponding to a neural network. If this neural network is the neural correlate of a quale[7], we ascribe the semantics of that quale to the Platonic set $\mathcal{K}_{d_\lambda p} d_\lambda p$, itself located in a virtual space. This quale is positioned in the rightmost box in Fig. 5.2.

## 6.3 Evolution

The characterization of the dynamics of memes (set pictures) and themata (the sets that are pictured) as an adaptive process, employing such Darwinian concepts as competition, selection, reproduction as well as fitness and genomics is also the subject of ongoing work (Miranker 2008) that finds motivation in the foundations developed here. We expect that variations of our development will provide mathematical foundations for the study of evolution driven by so-called *selfish replicators*, both genetic and mimetic (Dawkins 1979, Blackmore 1999).

## 6.4 Other directions

1. Study of the bisimulation of graphs, a notion that characterizes when two memes share a thema. We expect this to lead to a mathematical theory of memes and themata.
2. Model theoretic foundations of Aczel theory dealing with the consistency of the AFA with the Z-F Axioms from which FA has been deleted.
3. Study of the algebra of the set theoretic operators generated by the fundamental operators appearing in Table 3.1. An example of such an operator is $\mathcal{I} \cup \mathcal{BK}$.
4. Classification of the consciousness operators $\mathcal{K}$ and the connection of such a classification to Gödel's Incompleteness Theorem.
5. Examples and applications of the M-Z equation.

6. Application of these foundations to the grammar of programming languages.

7. Study of the trajectories generated by iterating application of a consciousness operator

---

[7] A quale is the perception of a color, an aroma…, or the perception of a feeling, such as hunger, fear... The neural correlate of a quale is the neural circuitry in the brain that is active when the quale is perceived. Some attribute the location of the quale to this circuitry.

## Appendix: Axioms of Set Theory

We make explicit use of the following axioms of set theory.

**Existence:** $\exists z(z = z).$

**Extensionality:** $\forall z(z \in a \leftrightarrow z \in b) \rightarrow a = b.$

**Pairing:** $\exists z[a \in z\ \&\ b \in z].$

**Union:** $\exists z(\forall x \in a)(\forall y \in x)(y \in z).$

**Comprehension:** $\exists z \forall x[x \in z \leftrightarrow x \in a\ \&\ \varphi(x)].$

Here $\varphi$ can be any formula in which the variable $z$ does not occur free.

Except for the axiom of existence these axioms along with the Axioms of Infinity, Collection, Power Set and Choice can be found in Aczel (1988). We do not state the latter four axioms since we use them only implicitly. Note that Aczel uses the name Axiom of Separation for the Axiom of Comprehension. The FA is stated as follows.

**Axiom of Foundation:** $\exists x(x \in a) \rightarrow (\exists x \in a)(\forall y \in x)\neg(y \in a).$

The FA is not included in the original Z-F list. It was proposed by Von Neumann. We don't use the FA, and we replace it by the AFA stated as follows.

**Anti-Foundation Axiom:** *Every graph has a unique decoration.*

The AFA, due to Aczel, is central to our development.

## Glossary[8]

**Terminology**
Experience/primary experience...a set *x*/primary layer when there are layers of experience
Consciousness...$\mathcal{K}x$, where $\mathcal{K}$ is a consciousness operator. See Semantic Thesis in §3.2
Awareness....$\mathcal{K}x$, where $\mathcal{K}$ is a consciousness operator. See Semantic Thesis in §3.2
Graph....a collection of nodes with certain pairs of the nodes specified as edges
Directed graph....a graph in which the nodal pairs are ordered (edges are directed)
Pointed graph....a directed graph with a distinguished node, the point
Accessible pointed graph (apg)....a pointed graph, every node of which is reachable from
    the point by a chain of directed edges
Decoration....the unique assignment (specified by (2.2)) of sets to the nodes of an apg
Picture of a set....the pointed graph in whose decoration, the set corresponds to the point

---

[8] For convenience, some of the definitions listed here are abbreviated. In such cases more complete definitions are found in the text.

Labeled graph….a graph with an arbitrary assignment of sets (the labels) to the nodes
Labeled decoration….a labeling dependent decoration of a graph (specified by (4.0))
Histogram….replaces a collection by a set as the domain of a set valued function
M-Z equation….specifies the labeled decoration of a graph arising from neural networks
Hebb's rule….specifies the synaptic weight change in a model neuron
McCulloch-Pitts equation….specifies the binary valued output of a model neuron

**Set types**
Collection….a set as defined by Cantor
Naïve set….another name for a collection
Set….a primitive construct, the subject of the Z-F axioms
Bona fide set….a set, emphasizing its being specified as a primitive defined by Z-F
Pure set….a set whose elements are sets, whose elements of elements are sets…
Path….a sequence of nodes (finite or infinite) linked by directed edges
Well-founded picture….a graph whose paths are finite (in particular, one without loops)
Non-well-founded picture….a graph with an infinite path
Well-founded set…. a set whose picture is well-founded
Non-well-founded set….a set whose picture is non-well-founded
Normal set….a set that does not contain itself
Abnormal set….a set that contains itself
Platonic set….a not physical set, a not computable set, a set located in a virtual space

**Classes**
Class….a collection of sets with a common property
Proper class….a class that is not a set
$\mathcal{U}$….the class or universe of sets
$\mathcal{A}$….the class of abnormal sets
$\mathcal{N}$….the class of normal sets

**Fundamental Operators**
$\mathcal{E}$….elimination
$\mathcal{I}$….identity
$\mathcal{B}$….brace, singleton
$\mathcal{R}$….Russell
$\mathcal{T}$….anti-Russell
$\mathcal{D}$….duality operator
$\mathcal{C}$….a particular consciousness operator

**Types of Operators**
$O$….a generic operator
$\mathcal{K}$….a generic consciousness operator
$\mathcal{K}_A$….a special class of consciousness operators parameterized by a set $A$
$\mathcal{K}^{diag}$….diagonalization of the family of operators $\mathcal{K}_A$
Selectors….operators $O$ with the following property: $x \subseteq y \Rightarrow Ox = x \cap Oy$

# References

Aczel, P. 1988. Non-Well-Founded *Sets*. CSLI Publications.

Aleksander, I., Dunmall, B. 2003. Axioms and Tests for the Presence of Minimal Consciousness in Agents. In *Machine Consciousness*, O. Holland, ed., Imprint Acad.

Aristotle. 1961. (Ross, D. Ed). *De Anima*, Oxford at the Clarendon Press.

Barwise. J. in the forward to Aczel 1988.

Bernays, P. 1954. A System of Axiomatic Set Theory, VII. *J. Symbolic Logic* 1981-86.

Blackmore, S. 1999. *The Meme Machine,* Oxford Univ. Press.

Cantor, G. 1895. Beiträge zur Begründung der Transfiniten Mengenlehre, 1. *Mathematische Annalen* 46; 481-512.

Courant, R., Robbins, H. 1941. *What is Mathematics.*

Dawkins, R. 1976. *The Selfish Gene.* Oxford University.

Descartes, R. 1637. *Discours sur la Methode.*

Dym, H., Mckean, H. 1972. *Fourier Series and Integrals.* Academic Press.

Fränkel, A. 1922. Zu den Grundlagender Cantor-Zermeloschen Mengenlehre. *Mathematische Annalen 86*; 230-237.

Frege, G. 1893. *Grundsetze der Arithmetik, begriffsschriftlich abgeleitet,* Vol. 1. Jena, Volume 2 published in 1903.

Gödel, K. see Jech.

Haykin, S. 1999. *Neural Networks a Comprehensive Foundation.* Prentice Hall.

Hilbert, D. 1900, Address to the International Congress of Mathematicians in Paris. See also Kaplansky. 1977. *Hilbert's Problems.* University of Chicago.

Hofstadter, D. 1979. *Gödel Escher Bach; An Eternal Golden Braid.*

Hrbacek, K., 1999. *Introduction to Set Theory.* M. Dekker.

Jech, T. 2002. *Set Theory*. Springer.

Kanger, S. 1957. *Provability in Logic.* Univ. of Stockholm: Almquist and Wiksell. Stockholm Studies in Philosophy.

Miranker, W. 2008, Memes and their Themata. In preparation.

Miranker, W., Zuckerman, G. 2008. Applications of the M-Z Equation. In preparation.

Penrose, R. 1989. *The Emperor's New Mind.* Oxford Univ. Press.

Plato. 360 BCE. *The Republic.* Translated by B. Jowett.

Russell, B., Whitehead, A., 1910-1913. *Philosophiae Naturalis Principia Mathematica.*

Schrödinger, E. 1958. *Mind and Matter*. Cambridge Univ. Press.

Solms, M. 1995. Chromosomes on the Couch. *Psychoanalytic Psychotherapy*, 9; 107-20.

Von Neumann, 1925. Eine Axiomatisierung der Mengenlehre. Jurnal für Reine und Angewandte Mathematik 155; 219-240.

Whitehead, A. see Russell.

Zermelo, E. 1908. Untersuchung uber die Grundlagen der Mengenlehre, I. *Mathematische Annalen* 65; 261-281.

**Power set:**     $\exists z \forall x \left[ (\forall u \in x)(u \in a) \to x \in z \right].$

**Infinity:**     $\exists z \left[ (\exists x \in z) \forall y \neg (y \in x) \ \& \ (\forall x \in z)(\exists y \in z)(x \in y) \right].$

**Collection:**     $(\forall x \in a) \exists y \varphi \to \exists z (\forall x \in a)(\exists y \in z) \varphi.$

**Choice:**     $(\forall x \in a) \exists y (y \in x)$
                         $\& \ (\forall x_1 \in a)(\forall x_2 \in a) \left[ \exists y (y \in x_1 \ \& \ y \in x_2) \to x_1 = x_2 \right]$
                         $\to \exists z (\forall x \in a)(\forall y \in x)(\forall u \in x) \left[ u \in z \leftrightarrow u = y \right].$